\title{\vspace{-0.7cm}Small Complete Minors Above the Extremal Edge Density}

\author{Asaf Shapira
\thanks{School of Mathematics, Tel-Aviv University, Tel-Aviv, Israel 69978, and Schools of Mathematics and Computer Science, Georgia Institute of Technology, Atlanta, GA 30332. Email:
{\tt asafico@tau.ac.il}. Supported in part by NSF Grant DMS-0901355, ISF Grant 224/11 and a Marie-Curie CIG Grant 303320.} \and Benny Sudakov
\thanks{Department of Mathematics, UCLA, Los Angeles, CA 90095. Email: {\tt bsudakov@math.ucla.edu}. Research supported
in part by NSF grant DMS-1101185, by AFOSR MURI grant FA9550-10-1-0569 and by a USA-Israel BSF
grant.
}
}

\date{}

\documentclass [letterpaper,11pt]{article}
\usepackage{amsfonts}

\newtheorem{theo}{Theorem}

\newtheorem{lemma}{Lemma}[section]
\newtheorem{definition}[lemma]{Definition}

\newtheorem{remark}[lemma]{Remark}
\newtheorem{claim}[lemma]{Claim}

\newcommand{\qed}{\hspace*{\fill} \rule{7pt}{7pt}}

\newcommand{\ignore}[1]{}

\begin{document}
\maketitle

\begin{abstract}

A fundamental result of Mader from 1972 asserts that a graph of high average degree contains a highly connected subgraph with roughly
the same average degree. We prove a lemma showing that one can strengthen Mader's result by replacing the notion
of high connectivity by the notion of vertex expansion.

Another well known result in graph theory states that for every integer $t$ there is a smallest real $c(t)$, such that every $n$-vertex graph with $c(t) n$ edges contains a
$K_t$-minor. Fiorini, Joret, Theis and Wood asked if an $n$-vertex graph $G$ has $(c(t)+\epsilon)n$ edges
then $G$ contains a $K_t$-minor of order at most $C(\epsilon)\log n$. We use our extension of Mader's theorem to prove that such a graph $G$ must contain a $K_t$-minor of order at most
$C(\epsilon)\log n\log\log n$. Known constructions of graphs with high girth show that this result is tight up to the $\log\log n$ factor.

\end{abstract}


\section{Introduction}\label{sec:intro}

\subsection{Graph minors and the main result}

All graphs considered here are finite and have no loops or parallel edges.
The order of a graph is the number of its vertices.
A graph $H$ is a minor of a graph $G$ if $H$ can be obtained from $G$ by a sequence
of edge deletions, vertex deletions and edge contractions. In this case we say
that $G$ has an $H$-minor. Let $K_t$ denote the complete graph on $t$ vertices. Since we will be mainly
interested in $K_t$-minors it will be easier for us to use the following equivalent definition of a $K_t$-minor.
A graph $G$ has a $K_t$-minor if $G$ contains $t$ vertex disjoint connected subgraphs
$S_1,\ldots,S_t$ and ${t \choose 2}$ paths $(P_{i,j})_{1 \leq i < j \leq t}$, such that $P_{i,j}$ connects $S_i$ to $S_j$,
each path $P_{i,j}$ is disjoint from all sets $S_k$ with $k \neq i,j$, and the paths $P_{i,j}$ are
internally vertex disjoint, that is, $P_{i,j}$ can only intersect with $P_{i,j'}$ or $P_{i',j}$ at its endpoint vertices.
This $K_t$-minor is called topological if each subgraph $S_i$ consists of a single
vertex.

The notion of minor is undoubtedly one of the most well studied topics in graph theory.
A central result in this area states that a linear number of edges is enough to force
the appearance of a $K_t$-minor. Formally, for every integer $t \geq 3$ define
\begin{equation}\label{ct}
c(t)=\min\{c: d(G) \geq c \mbox{ implies that $G$ has a $K_t$-minor}\}\;,
\end{equation}
where $d(G)=|E(G)|/|V(G)|$. Mader \cite{M2} has shown that the displayed set does indeed have a minimum
(that is, its infimum is a member of the set) and that $c(t) \leq 2^{t-3}$.
He later \cite{M3} obtained the improved bound $c(t) \leq 16t \log t$ (all logarithms in this paper are base 2). His result was improved by
Kostochka \cite{K} and Thomason \cite{Thomason1} who proved that $c(t)=\Theta(t\sqrt{\log t})$.
Thomason \cite{Thomason2} later proved an even sharper bound, showing that $c(t)=(\alpha+o(1))t\sqrt{\log t}$ for some explicit
constant $\alpha$, with the $o(1)$ term going to $0$ as $t \rightarrow \infty$.

Fiorini, Joret, Theis and Wood \cite{FJTW} raised the following problem: how many edges suffice to guarantee that a graph contains
not only a $K_t$-minor, but one which has {\em few vertices}? Observe that graphs with logarithmic girth (which can be constructed by deleting short cycles from
random graph $G(n,p), p=c/n$ or explicitly see, e.g., \cite{LPS})
show that for any constant $C$, there is a graph $G$ with $d(G)\geq C $ and no $K_3$ minor of order $o(\log n)$. We also need $d(G) \geq c(t) = (\alpha +o(1))t\sqrt{\log t}$ to guarantee
{\em some} $K_t$-minor. So the question boils down to finding the smallest constant $c> c(t)$, such that any graph $G$ with $d(G) \geq c$ contains a $K_t$-minor
of order $O(\log n)$. Fiorini et al. \cite{FJTW} proved
that if $d(G) \geq 2^{t-2}+\epsilon$ then $G$ has a $K_t$-minor of order $C(\epsilon)\log n$.
Note that the average degree here is exponentially larger than the one needed to guarantee a $K_t$-minor. This motivated Fiorini et al. \cite{FJTW} to ask if in fact any graph
$G$ with $d(G) \geq c(t)+\epsilon$ contains a $K_t$-minor of order $C(\epsilon)\log n$. That is, while
$c(t)n$ edges are sufficient (and necessary) to guarantee {\em some} $K_t$-minor, adding only $o(n)$ additional edges should force the appearance of the (asymptotically) smallest $K_t$-minor one can force even with $Cn$ edges, for
any constant $C$. Our main result in this paper comes very close to answering their question positively.

\begin{theo}\label{theomain} For every $\epsilon >0$ and integer $t \geq 3$ there exist $n_0=n_0(\epsilon,t)$, such that every $n$-vertex graph $G$
with $n\geq n_0$ and $d(G) \geq c(t)+\epsilon$ contains a $K_t$-minor of order $O(\frac{c(t)t^2}{\epsilon}\log n \log\log n)$.
\end{theo}

As we mentioned above, Thomason \cite{Thomason2} has shown that $c(t)=(\alpha+o(1))t\sqrt{\log t}$. The easier half of his result
shows that a random graph (on an appropriate number of vertices) with $d(G)=(1-\epsilon)c(t)$ has no $K_t$-minor.
Myers \cite{Myers} has later strengthened his result (see \cite{Myers} for the precise condition where this fact holds)
by showing that $n$-vertex graphs of density $p$ that do not contain a $K_t$-minor larger than the one we expect to find in a random graph $G(n,p)$ must be quasi-random.
So roughly speaking, random graphs are extremal with respect to the critical density where one expects to find a $K_t$-minor in {\em arbitrary} graphs.
A positive answer to the problem of Fiorini et al. \cite{FJTW} would thus show that random graphs are also
in some sense extremal with respect to the actual order of the $K_t$-minor we expect to find. Specifically, it is possible to show that the smallest $K_t$-minor in a random graph $G(n,p), p=c/n$ such that $d(G)=c/2 > c(t)$
has order $O(\log n)$.
The problem if \cite{FJTW} can thus be phrased as asking if in fact, appropriate random graphs maximize the order of the smallest $K_t$-minor among all graphs of density $c(t)+\epsilon$.

Let us finally mention an old conjecture of Erd\H{o}s, stating that a graph with $n^{1+\epsilon}$ edges contains a non-planar subgraph of size $C(\epsilon)$.
This conjecture was confirmed (in a very strong sense) by Kostochka and Pyber
\cite{KP} who proved that any graph with $4^{t^2}n^{1+\epsilon}$ edges contains a topological $K_t$-minor of size
$O(t^2\log t/\epsilon)$. So the problem of \cite{FJTW} that we study here is in some sense a strengthening of
the conjecture of Erd\H{o}s for $\epsilon=1/\log n$. Furthermore, as noted in \cite{FJTW}, one can adapt the argument of \cite{KP} to show that a graph $G$ with $d(G) \geq 16^t$ contains a $K_t$-minor of size $O(\log n)$.

\subsection{Expansion in graphs, proof overview and the key lemma}

We believe that an important aspect of this paper is the proof technique we employ here which relies on the notion of expansion in graphs
and might be applicable in other settings. Perhaps a good perspective on our approach comes from dense graphs, that is graphs with $cn^2$ edges. Probably the most powerful tool one has at his disposal when studying dense graphs is Szemer\'edi's regularity lemma \cite{Sz}, which asserts that {\em any} dense graph can be approximated by a graph consisting of a bounded number of quasi-random graphs. Since quasi-random graphs are much easier to work with, this lemma allows one to reduce a problem on arbitrary graphs to the same problem on quasi-random graphs. We refer the reader to \cite{KSim} for more details on the regularity lemma and its applications.

When it comes to sparse graphs, there is no analogue of the regularity lemma. But in recent years, a parallel paradigm has emerged, the underlying idea of which can be thought of as
stating that {\em any} sparse graph is close to being the disjoint union of expander graphs. While the regularity lemma supplies one notion of approximation/quasi-randomness for all
applications involving dense graphs, it seems like for sparse graphs different applications call for different notions of approximation and expansion. We refer the reader to
\cite{ABS,GR,Trevisan} for some examples where this paradigm was applied.

Just like graph minors, expansion is one of the most well studied topics in graph theory, with a remarkable number of applications in diverse areas such as theoretical computer
science, additive number theory and information theory (just to name a few). We will thus refrain from giving a detailed account and instead refer the reader to the surveys
\cite{HLW,KrivSu} for more details. There are several known results connecting expansion and existence of $K_t$-minors in graphs, see, e.g.,
\cite{AST, PRS, KR, KrivSu1}. In all these papers the goal was to maximize the value of $t$. Our task here is quite different, we want to minimize the number of vertices in
the minor, keeping $t$ fixed.

For the proof of Theorem \ref{theomain} we will need a very strong notion of expansion. The price will be that we will ask for a very weak notion of
approximation\footnote{In fact, what we will ask for in Lemma \ref{expansion} is just one subgraph with very good expansion properties, so this can hardly be called an approximation.
In Section \ref{sec:concluding} we will suggest a possible strengthening of Lemma \ref{expansion}, involving a stronger notion of approximation, and the possible applications of such a
lemma.}, which will turn out to be sufficient for proving Theorem \ref{theomain}. In what follows, for a set of vertices $S$ we use $N(S)$ to denote the neighborhood of $S$, that is
the set of vertices not in $S$ that are connected to at least one vertex in $S$. The notion of expansion we will use is the following:

\begin{definition}[$\delta$-Expander]\label{def:expander} An $m$-vertex graph $H$
is said to be a {\em $\delta$-expander} if for every integer $0 \leq d \leq \log\log m-1$ and $S \subseteq V(H)$ of order $|S| \leq m/2^{2^{d}}$ we have
\begin{equation}\label{eq:expansion}
|N(S)| \geq \frac{\delta 2^d}{\log m(\log\log m)^2}|S|
\end{equation}
\end{definition}

Observe that (disregarding the $(\log\log m)^2$ term) if $G$ is an $m$-vertex $\delta$-expander then sets of vertices of size $cm$ have vertex expansion about $1/\log m$ while sets of vertices of size $m^c$ have vertex expansion
$\Theta(1)$. The following lemma shows that we can indeed find a $\delta$-expander in any graph with sufficiently many edges.

\begin{lemma}[Key Lemma]\label{expansion}
If $G$ satisfies $d(G)=c$, then for every $0 < \delta \leq \frac{1}{256}$ we can find in $G$ a subgraph $H$, such that
$d(H) \geq (1-\delta)c$ and $H$ is a $\delta$-expander.
\end{lemma}

\begin{remark} Note that the only lower bound on the order of $H$ supplied by the lemma is $2(1-\delta)c$ which follows from the fact that
$|V(H)| \geq 2d(H) \geq 2(1-\delta)c$. Up to the $\delta$ error this is all that one can hope for since the graph might be a disjoint union of cliques
of order $2c$.
\end{remark}

In a nutshell, the proof of Theorem \ref{theomain} proceeds by first invoking Lemma \ref{expansion} on the graph $G$ thus obtaining a graph $H$ satisfying the expansion
properties of Definition \ref{def:expander}. We then show (see Lemma \ref{findminor1}) how one can find a small $K_t$-minor inside $H$, a task which is much easier given the fact that
$H$ has strong expansion properties. As we noted above, one can come up with different notions of expansion when studying sparse graphs. And indeed, in order to prove Theorem
\ref{theomain} we will actually have to prove another variant of Lemma \ref{expansion}, which uses a slightly different notion of expansion than the one defined in
(\ref{eq:expansion}). It might very well be possible to prove other variants of Lemma \ref{expansion}, suitable for tackling other problems. See the concluding remarks for another
possible variant which might have interesting applications.

As we mentioned in the abstract, Lemma \ref{expansion} can be thought of as a strengthening of Mader's Theorem (see \cite{D}, Theorem 1.4.3).
Indeed, Mader's Theorem states that any graph $G$ with $d(G) \geq 2k$ has a $k$-connected subgraph $H$ satisfying $d(H) \geq k$.
So Lemma \ref{expansion} gives a similar conclusion only it replaces the notion of $k$-connectivity with the stronger notion of vertex expansion.

\subsection{Organization} The rest of the paper is organized as follows. In Section \ref{sec:lemma} we prove Lemma \ref{expansion}. As we mentioned earlier, the proof can be easily
adapted to other settings, and indeed in Section \ref{sec:betterbound} we will prove another variant of this lemma. In Section \ref{sec:weaker} we prove a weaker version of Theorem
\ref{theomain}, giving a bound of order $O(\log n(\log\log n)^3)$. Then, in Section \ref{sec:betterbound}, we ``boost'' this weaker bound by removing a $(\log\log n)^2$ factor thus
obtaining the bound stated in Theorem \ref{theomain}. Section \ref{sec:concluding} contains some concluding remarks and open problems.

\section{Proof of Key Lemma}\label{sec:lemma}

In this section we prove Lemma \ref{expansion}. Given a graph $G$ and a subset of its vertices $U$, $G[U]$ denotes the subgraph of $G$ induced by the set $U$.
We start with the following simple yet crucial observation.

\begin{claim}\label{observation} If $G$ is an $n$-vertex graph satisfying $d(G)=c$, and $S \subset V(G)$ is such that $|N(S)| < \gamma |S|$, then either
$d(G[V \setminus S]) \geq c$ or $d(G[S \cup N(S)]) \geq (1-\gamma)c$.
\end{claim}

\paragraph{Proof:} Indeed, if $d(G[V \setminus S]) < c$,
$d(G[S \cup N(S)]) < (1-\gamma)c$ and $|N(S)| < \gamma|S|$ then we can bound
the number of edges in $G$ as follows,
$$
e(G) < c(n-|S|)+(1-\gamma)c(|S|+|N(S)|) \leq cn - \gamma c|S|+c|N(S)| \leq cn\;,
$$
contradicting the assumption that $G$ has $cn$ edges. $\qed$

\bigskip

In what follows, we say that a graph $G=(V,E)$ fails to be a $\delta$-expander at scale $d$ if there is a vertex set $S \subseteq V$ of size at most $m/2^{2^d}$ violating (\ref{eq:expansion}).

\paragraph{Proof of Lemma \ref{expansion}:} Set $G_0=G$ and
consider the following iterative process; if either $|G_t| \leq 256$ or if $G_t$ is a $\delta$-expander the process stops. Otherwise,
set $m=|G_t|$ and $d_t=d(G_t)$. Since $G_t$ is not a $\delta$-expander, there must be some $0 \leq d \leq \log\log m-1$ and a set $S_t \subseteq V(G_t)$ of size at most $m/2^{2^d}$
violating (\ref{eq:expansion}). We can deduce via Claim \ref{observation} (with $c=d_t$ and $\gamma=\frac{\delta 2^d}{\log
m(\log\log m)^2}$) that in this case either $d(G[V(G_t) \setminus S_t]) \geq d_t$ or $d(G[S_t \cup N(S_t)])  \geq (1-\frac{\delta 2^d}{\log
m(\log\log m)^2})d_t$. In the former case we set $G_{t+1}=G[V(G_t) \setminus S_t]$
and in the latter we set $G_{t+1}=G[S_t \cup N(S_t)]$. For brevity, in what follows we call the above
two cases, Case $1$ and Case $2$. Let $G'$ be the graph returned at the end of the process. We shall show that $d(G') \geq (1-\delta)c$, but we first
note that this will finish the proof. Indeed, if $|G'| > 256$ then by definition of the process we get that $G'$ must be a $\delta$-expander,
so $G'$ satisfies the requirements of the lemma.
If $|G'| \leq 256$ then one of its connected components, call it $G''$, must also satisfy $d(G'') \geq (1-\delta)c$. Since $\delta \leq 1/256$ and $G''$ is connected, we get
that $G''$ is (trivially) a $\delta$-expander, so we can return $G''$.

We now show that $d(G') \geq (1-\delta)c$. Recall that for each of the graphs $G_t$ in the process $d_t=d(G_t)$ (so $d_0=c$).
Note that at each iteration either $d_{t+1} \geq d_t$ (in Case $1$) or $d_{t+1}=\gamma_td_t$ for some $\gamma_t >0$.
Fix some $4 \leq k \leq \log\log n$ and suppose $r<r'$ are such that $G_r$ is the first graph in the process whose order is in the range $[2^{2^{k-1}},2^{2^k}]$ and $G_{r'}$ is the
last
such graph in the process. We wish to compute a lower bound on
$d_{r'+1}/d_r$, that is, a lower bound on the fraction of edge loss that can occur when ``passing'' through the interval $[2^{2^{k-1}},2^{2^k}]$.
Observe that if at some iteration $r \leq t \leq r'$ of the process, the graph $G_t$ fails to be a $\delta$-expander at scale $d$, then either $\gamma_t \geq 1$ (Case $1$), or $|G_{t+1}| \leq |G_t|/2^{2^d}$. For each $0\leq d \leq k$ let $a_d$ be the number of times Case $2$ occurred due to $G_t$ failing to be a $\delta$-expander at scale $d$. Then we have
$$
1 \leq  |G_{r'+1}| \leq |G_r|/\prod^{k}_{d=0}2^{a_d2^{d}}\leq 2^{2^k}/\prod^{k}_{d=0}2^{a_d2^{d}}\;,
$$
implying that
\begin{equation}\label{eq:size}
\sum^{k}_{d=0}a_d2^{d} \leq 2^{k}\;.
\end{equation}
As noted in the first paragraph of the proof, if at some iteration $r \leq t \leq r'$ Case $2$ happened due to the fact that $G_t$ failed to be a $\delta$-expander at scale $d$, then
we know that $\gamma_t \geq 1-\delta2^d/\log m(\log\log m)^2$ (recall that $m$ denotes the order of $G_t$). Since we are considering the case where $2^{2^{k-1}} \leq m \leq 2^{2^{k}}$ this means that in this range we have
\begin{equation}\label{gammat}
\gamma_t \geq 1-\delta2^d/2^{k-1}(k-1)^2\;.
\end{equation}
Using that $a_d$ was the number of times $\gamma_t$ satisfied (\ref{gammat}), we have that
\begin{eqnarray*}
d_{r'+1}/d_r=\prod^{r'}_{t=r}\gamma_t &\geq& \prod^{k}_{d=0}\left(1-\frac{\delta 2^{d}}{2^{k-1}(k-1)^2}\right)^{a_d} \geq
\left(1-\delta\sum^{k}_{d=0}\frac{a_d2^{d}}{2^{k-1}(k-1)^2}\right)\\ &\geq& \left(1-\frac{2\delta}{(k-1)^2}\right)\;,
\end{eqnarray*}
where the last inequality follows from (\ref{eq:size}).
We have thus established that $d_t$ goes down by a factor of at most $1-2\delta/(k-1)^2$ when the process passes through the interval $[2^{2^{k-1}},2^{2^k}]$.
Since we stop the process when $|G_t| \leq 256$, we only need to consider $k \geq 4$.
Hence, altogether $d_t$ can decrease by a factor of at most

\begin{equation}\label{finaleq}
\prod^{\log\log n}_{k=4}\left(1-\frac{2\delta}{(k-1)^2}\right) \geq 1-2\delta \sum^{\log\log n}_{k=4}\frac{1}{(k-1)^2} \geq 1-\delta.
\end{equation}
So when the process ends we indeed obtain a $\delta$-expander $G_{t'}$ satisfying $d(G_{t'}) \geq c(1-\delta)$.
$\qed$

\section{A Weaker Bound}\label{sec:weaker}

As mentioned above, our goal in this section is to prove the following slightly weaker version of Theorem \ref{theomain}.

\begin{lemma}\label{lemmaweak}
The following holds for every $\epsilon >0$, integer $t \geq 3$ and large enough $n \geq n_1(\epsilon,t)$.
Every $n$-vertex graph satisfying $d(G) \geq c(t)+\epsilon$ has a $K_t$-minor of order $O\big(\frac{c(t)t^2}{\epsilon}\log n(\log\log n)^3\big)$.
\end{lemma}

The main lemma we will use to prove Lemma \ref{lemmaweak} is the following result, which shows that one can find small $K_t$-minors in
$\delta$-expanders.

\begin{lemma}\label{findminor1} The following holds for all $\delta >0$ integer $t$ and large enough $m \geq m_0(\delta,t)$.
If $H$ is a $\delta$-expander on $m$ vertices, then it has a $K_t$-minor of order $O\big(\frac{t^2}{\delta}\log m(\log\log m)^3\big)$.
\end{lemma}

Let us first show how can one derive Lemma \ref{lemmaweak} from Lemma \ref{findminor1}

\paragraph{Proof of Lemma \ref{lemmaweak}:} We can clearly assume that $\epsilon < c(t)$ (otherwise we can replace $\epsilon$ with (say) $1/2$, which is smaller than $c(t)$ for all $t
\geq 3$).
Let $G$ be an $n$-vertex graph satisfying $d(G) = c(t)+\epsilon$.
Applying Lemma \ref{expansion} to $G$ with $\delta=\epsilon/3c(t)$ we obtain a subgraph $H$ satisfying
$$
d(H) \geq (1-\delta)(c(t)+\epsilon) > c(t)\;.
$$
Hence, if $|V(H)| < m_0(\delta,t)$ we get from the definition of $c(t)$ in (\ref{ct}) that $H$ has a $K_t$-minor
of order at most\footnote{Here we use the assumption that $n \geq n_1(\epsilon)$. Specifically, we need to take $n_1(\epsilon)=2^{m_0(\epsilon/3c(t),t)}$.}
$m_0(\delta,t)=m_0(\epsilon/3c(t),t) \leq \log(n)$.
If $|V(H)| \geq m_0(\delta,t)$ then we get from the second assertion of Lemma \ref{expansion} that $H$ must be
an $\frac{\epsilon}{3c(t)}$-expander. Hence, we can apply
Lemma \ref{findminor1} to find in $H$ a $K_t$-minor of order $O\big(\frac{c(t)t^2}{\epsilon}\log n(\log\log n)^3\big)$.
$\qed$

\bigskip

We now turn to prove Lemma \ref{findminor1}. For a subset of vertices $U$ of a graph $G$ we denote by $B_k(U)$ the ball of radius $k$ around $U$, i.e., the set of all vertices of $G$
which can be reached by a path of length at most $k$ from some vertex in $U$. We start with the following simple observation:

\begin{claim}\label{clm:findpath}
Suppose $U$ and $V$ are two vertex sets in an $m$-vertex graph $G$ such that the following condition holds for every integer $0 \leq d \leq \log\log m-1$;
whenever $|B_k(U)| \leq m/2^{2^{d}}$ we have
\begin{equation}\label{eq:findpath}
|N(B_{k+1}(U))| \geq \frac{\delta2^d}{10\log m (\log\log m)^2}|B_k(U)|\;,
\end{equation}
and a similar condition holds with respect to $V$. Then there is a path connecting $U$ to $V$ of length at most $\frac{20}{\delta}\log m (\log\log m)^3$.
\end{claim}

\paragraph{Proof:} It is clearly enough to show that $|B_k(U)| > m/2$ for some $k \leq \frac{10}{\delta}\log m (\log\log m)^3$.
So suppose $m/2^{2^{d}} \leq |B_k(U)| \leq m/2^{2^{d-1}}$. Then by (\ref{eq:findpath}) we either have $|B_{k+t}(U)| > m/2^{2^{d-1}}$ or
\begin{equation}\label{growball}
m/2^{2^{d-1}} \geq |B_{k+t}(U)| \geq |B_{k}(U)|\left(1+\frac{\delta2^{d-1}}{10\log m(\log\log m)^2}\right)^t \geq m /2^{2^{d}-\frac{\delta t2^{d-1}}{10\log m(\log\log m)^2}}\;.
\end{equation}
It is easy to check that the RHS of (\ref{growball}) is larger than the LHS when $t \geq \frac{10}{\delta}\log m (\log\log m)^2$. In other words, within
at most $\frac{10}{\delta}\log m (\log\log m)^2$ steps, the neighborhood around $U$ jumps from size larger than $m/2^{2^{d}}$ to size larger than $m/2^{2^{d-1}}$. Since there are only $\log\log m$ intervals $[m/2^{2^{d}},m/2^{2^{d-1}}]$ we see that $|B_k(U)| > m/2$ for some $k \leq \frac{10}{\delta}\log m (\log\log m)^3$.
$\qed$

\begin{definition}[$\gamma$-Expanding Ball]\label{def:cluster} Let $v$ be a vertex in a graph $G$. We say that $B_k(v)$ is a {\em $\gamma$-expanding ball} if for every
$1 \leq i \leq k-1$ we have $|B_{i+1}(v)| \geq |B_{i}(v)|(1+\gamma)$. We also call $v$ the {\em center} of the ball and $k$ the {\em radius}.
A set of vertices is said to be a $\gamma$-expanding ball if it is equal to an expanding ball $B_k(v)$ for some center $v$ and radius $k$.
\end{definition}

\begin{claim}\label{clm:sets1} The following holds for all $\delta>0$ integer $t$ and large enough $m \geq m_0(\delta,t)$.
If $G$ is a $\delta$-expander on $m$ vertices then one of the following holds;
\begin{enumerate}
\item $G$ has $t$ vertices $v_1,\ldots,v_t$ each of degree at least $\log^4m$.
\item Set $\gamma=\delta/5(\log\log m)^2$. Then $G$ contains $t$ disjoint vertex sets $S_1,\ldots,S_t$ such that $G[S_i]$ is a $\gamma$-expanding ball, $m^{1/5} \leq |S_i| \leq  m^{1/4}$ and every vertex in $G[S_i]$ has degree at most $\log^4m$.
\end{enumerate}
\end{claim}

\paragraph{Proof:} If $G$ has $t$ vertices $v_1,\ldots,v_t$ each of degree larger than $\log^4m$ then we have the first case of the lemma.
So suppose for the rest of the proof that $G$ has at most $t$ vertices of degree larger than $\log^4 m$. We need to show that we can pick sets $S_1,\ldots,S_t$ satisfying the second condition of the lemma. Let $T$ be the set containing the vertices of degree larger than $\log^4m$. Let us also say
that a set $S$ is {\em nice} if $G[S]$ is a $\gamma$-expanding ball and $m^{1/5} \leq |S| \leq  m^{1/4}$. Note that a nice set in $G \setminus T$ satisfies the second
condition of the lemma.

It is clearly enough to show that for any set $W \supseteq T$ of at most $m^{1/3}$ vertices, we can find in $G \setminus W$ a nice set $S$.
Once we know this, we can simply iteratively pick the sets $S_i$ one after the other where at iteration $i$ we will pick
$S_i$ from $G \setminus W_i$ with $W_i=(\bigcup_{j < i} S_j) \cup T$. Since $|T|,|S_1|,\ldots,|S_{i-1}| \ll m^{1/3}$ this set has size smaller than $m^{1/3}$.

So suppose to the contrary that there is a $W \supseteq T$ of size $m^{1/3}$ such that $G \setminus W$ has no nice set $S$ and set $G_0 = G \setminus W$.
As we assume that $W \supseteq T$, all the vertices in $G_0$ have degree at most $\log^4 m$. So pick a vertex $v_1 \in G_0$ and let $B_k(v_1)$ be the vertices at distance at most $k$
from $v_1$ in $G_0$. Let $k_1$ be the smallest integer such that $|B_{k_1+1}(v)| < |B_{k_1}(v)|(1+\gamma)$. We claim that $|B_{k_1+1}(v)| < m^{1/4}$. To see this observe that the fact
that all the vertices in $G_0$ have degree at most $\log^4 m$ implies that for every $k \leq k_1$ we have $|B_{k-1}(v)| \geq  |B_{k}(v)|/\log^4 m$ so if $|B_{k_1+1}(v)| \geq m^{1/4}$
then there must be a $k'_1$, such that $m^{1/5} \leq |B_{k'_1}(v)| \leq m^{1/4}$ and $|B_{i+1}(v)| \geq |B_{i}(v)|(1+\gamma)$ for all $1 \leq i \leq k'_1$. But in this case
$B_{k'_1}(v)$ would be a nice set. So setting $T_1 = B_{k_1}(v_1)$ we have $|T_1| < m^{1/4}$ and $|N_{G_0}(T_1)| < \gamma|T_1|$.

Let $G_1=G_0 \setminus T_1$. Take now another vertex $v_2 \in G_1$ and repeat the above process. We will eventually end up with a set
$T_2=B_{k_2}(v_2)$ of size smaller than $m^{1/4}$ satisfying $|N_{G_1}(T_2)| < \gamma|T_2|$. Define $G_2=G_1\setminus T_2$. We continue picking sets $T_i \subset G_{i-1}$ until the
first time
$|\bigcup_i T_i| > \sqrt{m}/2$. Since each set $T_i$ is of size at most $m^{1/4} \leq \sqrt{m}/2$ this means that
$\sqrt{m}/2 \leq |\bigcup_i T_i| \leq  \sqrt{m}$. Now, the fact that each of the sets $T_i$ satisfies $|N_{G_{i-1}}(T_i)| < \gamma|T_i|$ along with the facts that
$G_i=G_{i-1}\setminus T_i$ and $|\bigcup_i T_i| \leq
\sqrt{m}$ implies that
\begin{equation}\label{eq:exp2}
\left|N_{G_0}\left(\bigcup_i T_i\right)\right| \leq \sum_i|N_{G_{i-1}}(T_i)| < \sum_i \gamma|T_i|\leq \gamma\sqrt{m}=\delta\sqrt{m}/5(\log\log m)^2\;.
\end{equation}
Now recall that $G$ was assumed to be a $\delta$-expander and since $\sqrt{m}/2 \leq |\cup T_i| \leq  \sqrt{m}$ we must have (using $d=\log\log m -1$ in (\ref{eq:expansion}))
\begin{equation}\label{eq:exp1}
\left|N_G\left(\bigcup_i T_i\right)\right| \geq \delta\left|\bigcup_i T_i\right|/(\log\log m)^2 \geq  \delta\sqrt{m}/4(\log\log m)^2\;.
\end{equation}
Let us now recall that $G_0$ was obtained from $G$ by removing a set $W$ of no more than $m^{1/3}$ vertices. Hence, recalling (\ref{eq:exp2}) we see
$$
\left|N_G\left(\bigcup_i T_i\right)\right| \leq  |W| +\left|N_{G_0}\left(\bigcup_i T_i\right)\right| \leq m^{1/3} +\delta \sqrt{m}/5(\log\log m)^2 < \delta\sqrt{m}/4(\log\log m)^2\;,
$$
which contradicts (\ref{eq:exp1}).
$\qed$

\begin{claim}\label{clm:case1}
The following holds for all $\delta >0$ integer $t$ and large enough $m \geq m_0(\delta,t)$.
If $G$ is a $\delta$-expander of order $m$ that has $t$ vertices as in the first case of Claim \ref{clm:sets1} then it has
a $K_t$-minor of order $O\big(\frac{t^2}{\delta}\log m (\log\log m)^3\big)$.
\end{claim}

\paragraph{Proof:} Let $v_1,\ldots,v_t$ be the vertices satisfying the first assertion of Claim \ref{clm:sets1}. We will show that
in this case $G$ has a topological
$K_t$-minor of the required order. To do this, we show that we can find ${t \choose 2}$ paths $P_{i,j}$, where each path connects $v_i$ to $v_j$, has length at
most $\frac{20}{\delta}\log m (\log\log m)^3$ and is internally disjoint
from all other paths. Our plan is to successively find these paths by invoking Claim \ref{clm:findpath}. All we need to do is show that for any $v_i,v_j$ and any set $W$ of size at most $\log^2m$ the sets $\{v_i\},\{v_j\}$ satisfy (\ref{eq:findpath}) with respect to the
graph $G \setminus W$. Once we establish this fact, we will be able successively pick the paths $P_{i,j}$ via Claim \ref{clm:findpath} where at each iteration we will take $W$ to consist of the vertices $v_t$, $t \neq i,j$, together with the internal vertices of the paths $P_{i',j'}$ we have already picked.

We turn to show that $\{v_i\}$ satisfies (\ref{eq:findpath}) in the graph $G'=G \setminus W$ (the proof for $\{v_j\}$ is identical).
We first recall the assumption of the lemma that $v_i$ has degree $\log^4m$ in $G$. Since we assume that $|W| \leq \log^2 m$ we clearly have
$$
N_{G'}(v_i) \geq N_G(v_i) -|W| \geq  \log^4m - \log^2m \geq \frac12\log^4m\;.
$$
Hence, for all $k \geq 1$ we have
\begin{equation}\label{minsize}
|B_k(v_i)| \geq |B_1(v_i)| \geq \frac12\log^4m\;.
\end{equation}
Suppose now that $|B_k(v_i)| \leq m/2^{2^d}$. Since we assume that $G$ is a $\delta$-expander we deduce that
\begin{eqnarray*}
|N_{G'}(B_{k}(v_i))| &\geq& |N_{G}(B_{k}(v_i))| -|W| \geq \frac{\delta2^d }{\log m(\log\log m)^2}|B_{k}(v_i)| -\log^2m\\
&\geq& \frac{\delta2^d }{2\log m(\log\log m)^2}|B_{k}(v_i)|\;,
\end{eqnarray*}
where the last inequality follows from (\ref{minsize}). We thus get that $G'$ satisfies (\ref{eq:findpath}) with $U=\{v_i\}$ and $V=\{v_j\}$ so we can indeed find a path connecting $v_i$ to $v_j$ of length at most
$\frac{20}{\delta}\log m (\log\log m)^3$. $\qed$

\begin{claim}\label{clm:case2}
The following holds for all $\delta >0$ integer $t$ and large enough $m \geq m_0(\delta,t)$.
If $G$ is a $\delta$-expander of order $m$ that has $t$ sets as in the second case of Claim \ref{clm:sets1} then it has
a $K_t$-minor of order $O\big(\frac{t^2}{\delta}\log m (\log\log m)^3\big)$.
\end{claim}

\paragraph{Proof:} Let $S_1,\ldots,S_t$ be the sets satisfying the second condition of Claim \ref{clm:sets1}.
Recall that in this case each set $S_i$ is a $\gamma$-expanding ball around a vertex $v_i$ in the induced subgraph $G[S_i]$.
Also recall that $m^{1/5} \leq |S_i| \leq m^{1/4}$
and that all the vertices of $G[S_i]$ have degree at most $\log^4 m$.
This means that there is some $k_i$, such that the ball around $v_i$ in the graph $G[S_i]$ satisfies $\log^{4} m \leq |B_{k_i}(v_i)| \leq \log^8 m$.
For every $1 \leq i \leq t$ set $C_i=B_{k_i}(v_i)$ and note that $C_i$ is also a $\gamma$-expanding ball in $G[S_i]$ around the same center $v_i$.

We will shortly show that one can find ${t \choose 2}$ internally vertex disjoint paths $P_{i,j}$, where each $P_{i,j}$ connects $C_i$ to $C_j$, has length at most $\frac{20}{\delta}\log m (\log\log m)^3$ and avoids all the sets
$C_\ell, \ell \not =i,j$.
But let us first observe why this will conclude the proof. For each $1 \leq i \neq j \leq t$ let $Q_{i,j}$ be some path connecting
$v_i$ to the unique vertex of $P_{i,j}$ that belongs to $C_i$. For each $1 \leq i \leq t$ set $V_{i}=\cup_{j}Q_{i,j}$. Then each set $V_i$ is connected in $G$
and the paths $P_{i,j}$ are internally vertex disjoint and avoid the sets $V_\ell, \ell \not =i,j$ so contracting the sets $V_i$ indeed gives us a $K_t$-minor in $G$. As to the order of this minor, note that since
$C_i$ is a $\gamma$-expanding-ball (with $\gamma=\delta/5(\log\log m)^3$) of size at most $\log^8 m$ its radius is bounded by
$$
\log_{1+\gamma} |C_i| \leq \frac{2\log(\log^8 m)}{\gamma} \leq 80(\log\log m)^4/\delta\;,
$$
so each of the paths $Q_{i,j}$ is of length $o(\log m)$.
Hence the total size of the sets $V_i$ is much smaller than $\log m$. Thus, once we show how to find the above
mentioned paths $P_{i,j}$ with lengths bounded by $\frac{20}{\delta}\log m(\log\log m)^3$ we obtain a $K_t$-minor of order
at most
$$\log m+\frac{20t^2}{\delta}\log m (\log\log m)^3=O\Big(\frac{t^2}{\delta}\log m (\log\log m)^3\Big)\,.$$

Our plan is to show that one can simply iteratively pick the paths $P_{i,j}$ by successive applications of Claim \ref{clm:findpath}.
To this end, we need to show that after picking some of the paths, we can still find in the remaining graph another path. Since each of the paths
is of length at most $\frac{20}{\delta}\log m(\log\log m)^2 \ll \log^2m$, it is enough to show that for any set of vertices $W$ of size at most $\log^2m$, that is disjoint from each of
the sets $C_1,\ldots,C_t$, and such that for any $i \neq j$, the sets $C_i$ and $C_j$ satisfy (\ref{eq:findpath}) with respect to the graph
\begin{equation}\label{newG}
G'=G \setminus (W \cup \bigcup_{\ell \neq i,j} C_{\ell})\;.
\end{equation}
Once we know this, we will iteratively pick the paths $P_{i,j}$, where at each iteration we will take $W$ to be the union of the internal vertices of the paths we have already picked. To establish (\ref{eq:findpath}) we will need to consider two ``regimes'' of growth; the first is when $C_i$ grows within the set $S_i$, and the second, when it grows out of $S_i$.

Let us first make a simple (but crucial) observation about the sets $C_i$. Recall that $S_i=B_{k'_i}(v_i)$ in $G[S_i]$ for some vertex $v_i$ and $k'_i$,
that $|S_i| \geq m^{1/5}$ and that each vertex in $S_i$ has degree at most $\log^4 m$. This means that $k'_i \geq \log m/20\log\log m$. We also recall that $|C_i| \leq \log^8 m$ and
since $S_i$ is a $\gamma$-expanding ball with $\gamma \geq \delta/5(\log\log m)^3$ this means that in the induced subgraph
$G[S_i]$  we have that $C_i=B_{k_i}(v_i)$ for some $k_i \leq
\frac{80}{\delta}(\log\log m)^4$. Combining
the above facts about $k_i$ and $k'_i$ and setting $r_i=\log m/30\log\log m$ we have in graph $G[S_i]$
\begin{equation}\label{ballinball}
B_{r_i}(C_i) \subseteq S_i\;.
\end{equation}

We now turn to show that $C_i$ satisfies (\ref{eq:findpath}) with respect to $G'$ defined in (\ref{newG}).
Since $S_i$ is a $\gamma$-expanding ball we get from (\ref{ballinball}) that
$|N_{G[S_i]}(C_i)| \geq \delta|C_i|/5(\log\log m)^2$.
Recall now that $S_i$ is disjoint from $\bigcup_{r \neq i,j} C_r$, and that (\ref{ballinball}) implies that $N_{G[S_i]}(C_i) \subseteq S_i$. Hence
\begin{eqnarray*}
|N_{G'}(C_i)| &\geq&  |N_{G[S_i]}(C_i)|- |W| \geq \delta|C_i|/5(\log\log m)^2- \log^2m \\
&\geq& \delta|C_i|/10(\log\log m)^2\;,
\end{eqnarray*}
where the last inequality uses the fact that $|C_i| \geq \log^4 m$. In other words, we have in $G'$
\begin{equation}\label{expandinSi}
|B_1(C_i)|\geq \left(1+\frac{\delta}{10(\log\log m)^2}\right)|C_i|\;.
\end{equation}
Recalling (\ref{ballinball}) we can continue inductively and get that for all $k \leq r_i$ we have
\begin{eqnarray}
\label{beny1}
|N_{G'}(B_{k}(C_i))| &\geq& |N_{G[S_i]}(B_{k}(C_i))| - |W| \geq \delta|B_{k}(C_i)|/5(\log\log m)^2- \log^2m \nonumber\\
&\geq& \delta|B_{k}(C_i)|/10(\log\log m)^2\;,
\end{eqnarray}
where in the last inequality we use the fact that $|B_{k}(C_i)| \geq |C_i| \geq \log^4 m$. It follows that $C_i$ satisfies
(\ref{eq:findpath}) for all $k \leq r_i$ (in the graph $G'$). To see this,
note that by (\ref{ballinball}) $B_{k}(C_i) \subset B_{r_i}(C_i) \subseteq S_i$ implying that $|B_{k}(C_i)| \ll \sqrt{m}$. This means
that in (\ref{eq:findpath}) the relevant $d$ is $\log\log m -1$,
implying that we should show that $|N_{G'}(B_{k}(C_i))|\geq \delta|B_{k}(C_i)|/20(\log\log m)^2 $ as we indeed derived in (\ref{beny1}).
From (\ref{beny1}) we also have
\begin{equation}\label{largeball}
|B_{r_i}(C_i)| \geq |C_i|\left(1+\frac{\delta}{10(\log\log m)^2}\right)^{r_i} \geq \left(1+\frac{\delta}{10(\log\log m)^2}\right)^{\frac{\log m}{30\log\log m}} \geq \log^{10}m\;.
\end{equation}
Consider now some $k > r_i$ and suppose $|B_k(C_i)| \leq m/2^{2^{d}}$. Then since $G$ is assumed to be a $\delta$-expander we have
\begin{eqnarray*}
|N_{G'}(B_{k}(C_i))| &\geq& |N_{G}(B_{k}(C_i))| - |W|-\sum_{\ell}|C_{\ell}| \geq \delta2^d|B_{k}(C_i)|/(\log\log m)^2- \log^9m \\
&\geq& \delta2^d|B_{k}(C_i)|/2(\log\log m)^2\;,
\end{eqnarray*}
where the last inequality follows from (\ref{largeball}) which tells us that for any $k \geq r_i$ we have $|B_{k}(C_i)| \geq |B_{r_i}(C_i)| \geq \log^{10} m$.
So $C_i$ satisfies (\ref{eq:findpath}) for all $d$.
$\qed$

\bigskip

\paragraph{Proof of Lemma \ref{findminor1}:} Immediate from Claims \ref{clm:sets1}, \ref{clm:case1} and \ref{clm:case2}. $\qed$

\section{Proof of Theorem \ref{theomain}}\label{sec:betterbound}

As discussed in Section \ref{sec:intro}, one can come up with different variants of Lemma \ref{expansion} involving different notions
of expansion. And indeed, to prove Theorem \ref{theomain} we will need to redo most of Section \ref{sec:weaker} with a different notion of expansion, defined as follows.

\begin{definition}[$(\delta,n)$-Expander]\label{def:nexpander} An $m$-vertex graph $H$
is said to be a {\em $(\delta,n)$-expander} if for every integer $0 \leq d \leq \log\log m-1$ and $S \subseteq V(H)$ of size $|S| \leq m/2^{2^{d}}$ we have
\begin{equation}\label{eq:expansion2}
|N(S)| \geq \frac{\delta 2^d}{\log n}|S|
\end{equation}
\end{definition}

To apply the above notion of expansion we will need the following two lemmas, which are appropriate variants of Lemmas \ref{expansion} and \ref{findminor1}.

\begin{lemma}\label{lem:imp-bound} Let $G$ be an $n$-vertex graph satisfying $d(G)=c$. Then
for every $\delta >0$, the graph $G$ has a subgraph $G'$, such that $d(G') \geq (1-2\delta)c$ and $G'$ is a $(\delta,n)$-expander
\end{lemma}

\begin{lemma}\label{lem:findminor} The following holds for all $0 < \delta < 1/4$ integer $t$ and large enough $m \geq m_0(\delta,t)$.
If $G$ is an $m$-vertex $(\delta,n)$-expander with $2^{\log n/(\log\log n)^2} \leq m \leq n$, then $G$ contains
a $K_t$-minor with at most $O\big(\frac{c(t)t^2}{\epsilon}\log n \log\log n\big)$ vertices.
\end{lemma}

Let us first show how to derive Theorem \ref{theomain} from the above two lemmas and the result of the previous section.

\paragraph{Proof of Theorem \ref{theomain}:}
We can clearly assume that $\epsilon < c(t)$ (otherwise we can replace $\epsilon$ with (say) $1/2$, which is smaller than $c(t)$ for all $t \geq 3$).
Let $G$ be an $n$-vertex graph satisfying $d(G) = c(t)+\epsilon$.
Apply Lemma \ref{lem:imp-bound} on $G$ with $\delta=\epsilon/8c(t)$. Suppose first that the graph $G'$ returned by the lemma
is of order $m \leq 2^{\log n/(\log\log n)^2}$. Then the assertion of the lemma guarantees that $d(G') \geq (1-2\delta)(c(t)+\epsilon)\geq c(t)+\epsilon/2$.
If $m \leq n_1(\epsilon/2,t)$ we return $G'$. The definition of $c(t)$ then guarantees that $G'$ has a $K_t$-minor, whose order is
at most\footnote{Here we use the assumption that $n \geq n_0(\epsilon,t)$. Specifically, we take $n_0(\epsilon,t)=2^{n_1(\epsilon/2,t)}$, where $n_1(\epsilon,t)$ is the function used
in Lemma \ref{lemmaweak}.} $n_1(\epsilon/2,t) \leq \log n$ and we are done.
If $n_1(\epsilon/2,t) \leq  m \leq  2^{\log n/(\log\log n)^2}$ then Lemma \ref{lemmaweak} guarantees that $G'$ contains a $K_t$-minor of order at most
$O\big(\frac{c(t)t^2}{\epsilon}\log m
(\log\log m)^3\big)=O\big(\frac{c(t)t^2}{\epsilon}\log n \log\log n\big)$ as needed. Finally, if $m \geq 2^{\log n/(\log\log n)^2}$ then
the fact that $G'$ must be a $(\delta,n)$ expander guarantees, together with Lemma\footnote{Here we again use the assumption that $n \geq n_0(\epsilon,t)$.
Specifically, it is enough to take $n_0(\epsilon,t)=(m_0(\delta,t))^{\log m_0(\delta,t)}$, where $\delta=\epsilon/8c(t)$ and $m_0(\delta,t)$ is the function in
Lemma \ref{lem:findminor}} \ref{lem:findminor}, that $G'$ contains a $K_t$-minor of order $O\big(\frac{c(t)t^2}{\epsilon}\log n \log\log n\big)$ thus completing the proof.
$\qed$

\bigskip

Let us now turn to prove Lemmas \ref{lem:imp-bound} and \ref{lem:findminor}.
In what follows, we say that a graph $G=(V,E)$ fails to be a $(\delta,n)$-expander at scale $d$ if there is a vertex set $S \subseteq V$ of size at most $m/2^{2^d}$ violating (\ref{eq:expansion2}).

\paragraph{Proof of Lemma \ref{lem:imp-bound}:} Set $G_0=G$ and consider the following
iterative process; if either $|G_t| \leq 4$ or if $G_t$ is a $(\delta,n)$-expander the process returns $G_t$. Otherwise,
set $m=|G_t|$ and $d_t=d(G_t)$. Since $G_t$ is not a $(\delta,n)$-expander, there must be some $0 \leq d \leq \log\log m-1$ and a set $S_t \subseteq V(G_t)$ of size at most $m/2^{2^d}$
violating
(\ref{eq:expansion2}). We can deduce via Claim \ref{observation} (with $c=d_t$ and $\gamma=\delta 2^d/\log n$) that in this case either $d(G[V(G_t) \setminus S_t]) \geq d_t$
or $d(G[S_t \cup N(S_t)])  \geq (1-\delta 2^d/\log n)d_t$. In the former case we set $G_{t+1}=G[V(G_t) \setminus S_t]$
and in the latter we set $G_{t+1}=G[S_t \cup N(S_t)]$. For brevity, in what follows we call the above
two cases, Case $1$ and Case $2$. Let $G'$ be the graph returned at the end of the process.
We will now show that $d(G') \geq (1-2\delta)c$. Note that this assertion implies that if $|V(G')| > 4$ then the definition
of the process implies that $G'$ is a $(\delta,n)$-expander and if $|V(G')| \leq 4$, then since $\delta < 1/4$ one of $G'$ connected
components, say $G''$, is (trivially) a $\delta$-expander satisfying $d(G'') \geq (1-2\delta)c$.

Recall that for each of the graphs $G_t$ in the process $d_t=d(G_t)$ (so $d_0=c$). Note that at each iteration either $d_{t+1} \geq d_t$ (in Case $1$) or $d_{t+1}=\gamma_td_t$ for
some $\gamma_t >0$. Fix some $1 \leq k \leq \log\log n$ and suppose $r<r'$ are such that $G_r$ is the first graph in the
process whose size is in the range $[2^{2^{k-1}},2^{2^k}]$ and $G_{r'}$ is the last such graph in the process. We wish to compute a lower bound on
$d_{r'+1}/d_r$, that is, a lower bound on the fraction of edge loss that can occur when ``passing'' through the interval $[2^{2^{k-1}},2^{2^k}]$.  Observe that if at some iteration $r \leq t \leq r'$ of the process, the graph $G_t$ fails to be a $(\delta,n)$-expander at scale $d$, then either $\gamma_t \geq 1$ (Case $1$), or $|G_{t+1}| \leq |G_t|/2^{2^d}$. For each $0\leq d \leq k$ let $a_d$ be the number of times Case $2$ occurred due to $G_t$ failing to be a $(\delta,n)$-expander at scale $d$. Then we have
$$
1 \leq  |G_{r'+1}| \leq |G_r|/\prod^{k}_{d=0}2^{a_d2^{d}}\leq 2^{2^k}/\prod^{k}_{d=0}2^{a_d2^{d}}\;,
$$
implying that
\begin{equation}\label{eq:size1}
\sum^{k}_{d=0}a_d2^{d} \leq 2^{k}\;.
\end{equation}
Now, if at some iteration $r \leq t \leq r'$ Case $2$ happened due to the fact that $G_t$ failed to be a $(\delta,n)$-expander at scale $d$, then (as noted above) we know that in this case
\begin{equation}\label{gammat1}
\gamma_t \geq 1-\delta2^d/\log n\;.
\end{equation}
Using that $a_d$ was the number of times $\gamma_t$ satisfied (\ref{gammat1}), we have that
\begin{eqnarray*}
d_{r'+1}/d_r=\prod^{r'}_{t=r}\gamma_t &\geq& \prod^{k}_{d=0}\left(1-\frac{\delta 2^{d}}{\log n}\right)^{a_d}
\geq \left(1-\delta\sum^{k}_{d=0}\frac{a_d2^{d}}{\log n}\right)\\
&\geq& \left(1-\frac{\delta2^{k}}{\log n}\right)\;,
\end{eqnarray*}
where the last inequality follows from (\ref{eq:size1}).
We have thus established that $d_t$ goes down by a factor of at most $1-\delta2^{k}/\log n$ when the process passes through the interval $[2^{2^{k-1}},2^{2^k}]$. Hence, altogether $d_t$ can decrease by a factor of at most

\begin{equation}\label{finaleq1}
\prod^{\log\log n}_{k=1}\left(1-\frac{\delta2^{k}}{\log n}\right) \geq 1-\frac{\delta}{\log n} \sum^{\log\log n}_{k=1}2^k \geq 1-2\delta.
\end{equation}
So when the process ends we indeed obtain a $(\delta,n)$-expander $G_{t'}$ satisfying $d(G_{t'}) \geq c(1-2\delta)$.
$\qed$

\bigskip

As it turns out the proof of Lemma \ref{lem:findminor} can be obtained by repeating almost verbatim the proofs of Claims \ref{clm:findpath}, \ref{clm:sets1}, \ref{clm:case1} and \ref{clm:case2}, while replacing the notion of $\delta$-expander with $(\delta,n)$-expander, and making some minor adaptations to the calculations. Hence we only state the appropriate variants of these claims. The detailed proofs of these claims can be found in the appendix of the arxiv version of this paper.
The proof of Lemma \ref{lem:findminor} follows immediately from Claims \ref{clm:sets11}, \ref{clm:case11} and
\ref{clm:case21}.

\begin{claim}\label{clm:findpath1}
Suppose $U$ and $V$ are two vertex sets in an $m$-vertex graph $G$ such that the following condition holds for every integer $0 \leq d \leq \log\log m-1$;
whenever $|B_k(U)| \leq m/2^{2^{d}}$ we have
\begin{equation}\label{eq:findpath1}
|N(B_{k+1}(U))| \geq \frac{\delta2^d}{10\log n }|B_k(U)|\;,
\end{equation}
and a similar condition holds with respect to $V$. Then there is a path connecting $U$ to $V$ of length at most $\frac{20}{\delta}\log n (\log\log n)$.
\end{claim}

\begin{claim}\label{clm:sets11} The following holds for all $\delta >0$ integer $t$ and large enough $n \geq n_0(\delta,t)$.
If $G$ is a $(\delta,n)$-expander on $m$ vertices with $2^{\log n/(\log\log n)^2} \leq m \leq n$ then one of the following holds;
\begin{enumerate}
\item $G$ has $t$ vertices $v_1,\ldots,v_t$ each of degree at least $\log^4m$.
\item Set $\gamma=\frac{\delta\log m}{5\log n}$. Then $G$ contains $t$ disjoint vertex sets $S_1,\ldots,S_t$ such that $G[S_i]$ is a $\gamma$-expanding ball, $m^{1/5} \leq |S_i| \leq  m^{1/4}$ and every vertex in $G[S_i]$ has degree at most $\log^4m$.
\end{enumerate}
\end{claim}

\begin{claim}\label{clm:case11}
The following holds for all $\delta >0$ integer $t$ and large enough $n \geq n_0(\delta,t)$.
If $G$ is a $(\delta,n)$-expander on $2^{\log n/(\log\log n)^2} \leq m \leq n$ vertices which has $t$ vertices as in the first case of
Claim \ref{clm:sets11} then it
has a $K_t$-minor of order $O\big(\frac{t^2}{\delta}\log n \log\log n\big)$.
\end{claim}

\begin{claim}\label{clm:case21}
The following holds for all $\delta >0$ integer $t$ and large enough $n \geq n_0(\delta,t)$.
If $G$ is a $(\delta,n)$-expander on $2^{\log n/(\log\log n)^2} \leq m \leq n$ vertices which has $t$ sets as in the second case of Claim
\ref{clm:sets11} then it has a $K_t$-minor of order $O\big(\frac{t^2}{\delta}\log n \log\log n\big)$.
\end{claim}

\section{Concluding Remarks and Open Problems}\label{sec:concluding}

\begin{itemize}

\item Of course, the obvious open problem is to close the $O(\log\log n)$ gap between the upper bound of Theorem \ref{theomain} and the simple lower bound mentioned
in Section \ref{sec:intro}. It seems that
we have pushed our approach to the limit, and that in order to obtain a $K_t$-minor of order $C(\epsilon)\log n$ some new ideas
are needed.

We note that the proof of \cite{FJTW} showing that $(2^{t-2}+\epsilon)n$ edges suffice to guarantee a $K_t$-minor of order $C(\epsilon)\log n$
proceeds by applying the well known argument showing that $c(t) \leq 2^{t-3}$. In some sense,
we are able to prove our bound already for $d(G) \geq c(t)+\epsilon$ since we do {\em not} implicitly prove any bound on $c(t)$ but instead
rely on it as a black-box. Having said that, it might very well be possible
to use one of the proofs that give a tight bound on $c(t)$ to improve our lower bound.

\item A natural variant of the problem studied here concerns topological $K_t$-minors.
As noted in \cite{FJTW}, the method of \cite{KP} can probably be used to show that a graph
with $4^{t^2}n$ edges has a topological $K_t$-minor of order $O(\log n)$. However, it is known (see \cite{BT} and \cite{D})
that there is a smallest $c'(t)=\Theta(t^2)$ such that any graph with $c'(t)n$ edges has a topological $K_t$-minor.
It thus seems reasonable to conjecture that if $G$ is an $n$-vertex graph with $(c'(t)+\epsilon)n$ edges
then $G$ must contain a topological $K_t$-minor of order at most $C(\epsilon)\log n$.
We suspect that the ideas and tools in the present paper should allow
one to prove (at least) a weaker version of this conjecture, giving a topological $K_t$-minor of order at most $C(\epsilon)\log n\mbox{~poly}(\log\log n)$.

\item Let us say that a graph $G$ is $\epsilon$-far from being $K_t$-minor free if one should remove from $G$ at least $\epsilon n$
edges in order to turn $G$ into a $K_t$-minor free graph. We conjecture that if $G$ is $\epsilon$-far from being $K_t$-minor free then
$G$ has a $K_t$-minor of order $C(\epsilon)\log n$. At the moment, we are unable to prove even a weaker result in which the $K_t$-minor is of order
$C(\epsilon)\log n\mbox{~poly}(\log\log n)$. One approach to proving such a result would be to prove a strengthened version of Lemma \ref{expansion} stating
that for any $\delta >0$, and large enough graph $G$ one can remove from $G$ at most $\delta n$ edges and thus obtain a graph $G'$ in which every connected
component is a $\delta$-expander. Note that Lemma \ref{findminor1} together with the above variant of Lemma \ref{expansion} would immediately resolve the
above conjecture (with the slightly weaker bound $O(\log n\mbox{~poly}(\log\log n))$). Finally, observe that any graph with $(c(t)+\epsilon)n$ edges is by definition $\epsilon$-far from being $K_t$-minor free, hence
this conjecture strengthens the problem raised in \cite{FJTW}. One can of course raise the same conjecture with respect to topological minors.

\end{itemize}

\noindent
{\bf Acknowledgment.}\, We would like to thank Jie Ma for stimulating discussions.

\section{Missing Proofs from Section \ref{sec:betterbound}}

\paragraph{Proof of Claim \ref{clm:findpath1}:}
It is clearly enough to show that $|B_k(U)| > m/2$ for some $k \leq \frac{10}{\delta}\log n \log\log n$.
So suppose $m/2^{2^{d}} \leq |B_k(U)| \leq m/2^{2^{d-1}}$. Then by (\ref{eq:findpath1}) we either have $|B_{k+t}(U)| > m/2^{2^{d-1}}$ or
\begin{equation}\label{growball1}
m/2^{2^{d-1}} \geq |B_{k+t}(U)| \geq |B_{k}(U)|\left(1+\frac{\delta2^{d-1}}{10\log n}\right)^t \geq m /2^{2^{d}-\frac{\delta t2^{d-1}}{10\log n}}\;.
\end{equation}
It is easy to check that the RHS of (\ref{growball1}) is larger than the LHS when $t \geq \frac{10}{\delta}\log n $. In other words, within
at most $\frac{10}{\delta}\log n$ steps, the neighborhood around $U$ jumps from size larger than $m/2^{2^{d}}$ to size
larger than $m/2^{2^{d-1}}$. Since there are only $\log\log m \leq \log\log n$ intervals $[m/2^{2^{d}},m/2^{2^{d-1}}]$ we see that $|B_k(U)| > m/2$ for some $k \leq
\frac{10}{\delta}\log n \log\log n$. $\qed$

\medskip

\paragraph{Proof of Claim \ref{clm:sets11}:}
If $G$ has $t$ vertices $v_1,\ldots,v_t$ each of degree larger than $\log^4m$ then we have the first case of the lemma. So suppose for the rest of the proof that $G$
has at most $t$ vertices of degree larger than $\log^4 m$. We need to show that we can pick sets $S_1,\ldots,S_t$ satisfying the second condition of the lemma. Let $T$ be the set
containing the vertices of degree larger than $\log^4m$. Let us also say that a set $S$ is {\em nice} if $G[S]$ is a $\gamma$-expanding ball and $m^{1/5} \leq |S| \leq m^{1/4}$. Note that a
nice set in $G \setminus T$ satisfies the second condition of the lemma.

It is clearly enough to show that for any set $W \supseteq T$ of at most $m^{1/3}$ vertices, we can find in $G \setminus W$ a nice set $S$.
Once we know this, we can simply iteratively pick the sets $S_i$ one after the other where at iteration $i$ we will pick
$S_i$ from $G \setminus W_i$ with $W_i=(\bigcup_{j < i} S_j) \cup T$. Since $|T|,|S_1|,\ldots,|S_{i-1}| \ll m^{1/3}$ this set has size smaller than $m^{1/3}$.

So suppose to the contrary that there is a $W \supseteq T$ of size $m^{1/3}$ such that $G \setminus W$ has no nice set $S$ and let $G_0 = G \setminus W$. As we assume that $W \supseteq
T$, all the vertices in $G_0$ have degree at most $\log^4 m$. So pick a vertex $v_1 \in G_0$ and let $B_k(v_1)$ be the vertices at distance at most $k$ from $v_1$ in $G_0$. Let $k_1$
be
the smallest integer such that $|B_{k_1+1}(v)| < |B_{k_1}(v)|(1+\gamma)$. We claim that $|B_{k_1+1}(v)| < m^{1/4}$. To see this observe that the fact that all the vertices in $G_0$ have
degree at most $\log^4 m$ implies that for every $k \leq k_1$ we have $|B_{k-1}(v)| \geq |B_{k}(v)|/\log^4 m$ so if $|B_{k_1+1}(v)| \geq m^{1/4}$ then there must be a $k'_1$, such that
$m^{1/5} \leq |B_{k'_1}(v)| \leq m^{1/4}$ and $|B_{i+1}(v)| \geq |B_{i}(v)|(1+\gamma)$ for all $1 \leq i \leq k'_1$. But in this case $B_{k'_1}(v)$ would be a nice set. So
setting $T_1 = B_{k_1}(v_1)$ we have $|T_1| < m^{1/4}$ and $|N_{G_0}(T_1)| < \gamma|T_1|$.

Let $G_1=G_0\setminus T_1$. Take now another vertex $v_2 \in G_1$ and repeat the above process. We will eventually end up with a set
$T_2=B_{k_2}(v_2)$ of size smaller than $m^{1/4}$ satisfying $N_{G_1}|T_2| < \gamma|T_2|$. Let $G_2=G_1\setminus T_2$. We continue picking sets $T_i\subset G_{i-1}$ until the first
time
$|\bigcup_i T_i| > \sqrt{m}/2$. Since each set $T_i$ is of size at most $m^{1/4} \leq \sqrt{m}/2$ this means that $\sqrt{m}/2 \leq |\bigcup_i T_i| \leq  \sqrt{m}$. Now, the fact that
each of the sets $T_i$ satisfies $|N_{G_{i-1}}(T_i)| < \gamma|T_i|$ along with the facts that $G_i=G_{i-1}\setminus T_i$ and $|\bigcup_i T_i| \leq  \sqrt{m}$ implies that
\begin{equation}\label{eq:exp21}
\left|N_{G_0}\left(\bigcup_i T_i\right)\right| \leq \sum_i|N_{G_{i-1}}(T_i)| < \sum_i \gamma|T_i| \leq \gamma\sqrt{m}=\delta\sqrt{m}\log m/5\log n\;.
\end{equation}
Now recall that $G$ was assumed to be a $(\delta,n)$-expander and since $\sqrt{m}/2 \leq |\cup T_i| \leq  \sqrt{m}$ we must have (using $d=\log\log m -1$ in (\ref{eq:expansion2}))
\begin{equation}\label{eq:exp11}
\left|N_G\left(\bigcup_i T_i\right)\right| \geq \delta\left|\bigcup_i T_i\right|\frac{\log m}{2\log n} \geq \delta\sqrt{m}\log m/4\log n \;.
\end{equation}
Note that $G_0$ was obtained from $G$ by removing a set $W$ of no more than $m^{1/3}$ vertices and that by our assumption $m\gg \log^6 n$ . Hence, using (\ref{eq:exp21}) we have
$$
\left|N_G\left(\bigcup_i T_i\right)\right| \leq  |W| +\left|N_{G_0}\left(\bigcup_i T_i\right)\right| \leq m^{1/3} +\delta\sqrt{m}\log m/5\log n < \delta\sqrt{m}\log m/4\log n\;,
$$
which contradicts (\ref{eq:exp11}).
$\qed$

\medskip

\paragraph{Proof of Claim \ref{clm:case11}:}
Let $v_1,\ldots,v_t$ be the vertices satisfying the first assertion of Claim \ref{clm:sets1}. We will show that in this case $G$ has a
topological $K_t$-minor of the required order. To do this we show that we can find ${t \choose 2}$ paths $P_{i,j}$, where each path connects
$v_i$ to $v_j$, has length at most $\frac{20}{\delta}\log n \log\log n$ and is internally disjoint from all other paths. Our plan is to successively find
these paths by invoking Claim \ref{clm:findpath1}. All we need to do is show that for any $v_i,v_j$ and any set $W$ of size at most $\log^2 n$ the sets $\{v_i\},\{v_j\}$ satisfy
(\ref{eq:findpath1}) with respect to the graph $G \setminus W$. Once we establish this fact, we will be able successively pick the paths $P_{i,j}$ via Claim \ref{clm:findpath1} where
at each iteration we will take $W$ to consist of the vertices $v_t$, $t \neq i,j$, together with the internal vertices of the paths $P_{i',j'}$ we have already picked.

We turn to show that $\{v_i\}$ satisfies (\ref{eq:findpath1}) in the graph $G'=G \setminus W$ (the proof for $\{v_j\}$ is identical).
We first recall the assumption of the lemma that $v_i$ has degree $\log^4m$ in $G$. Since we assume that $|W| \leq \log^2 n$ and $\log m \geq \log^{1-o(1)} n$ we have
$$
N_{G'}(v_i) \geq N_G(v_i) -|W| \geq  \log^4m - \log^2n \geq \frac12\log^4m\;.
$$
Hence, for all $k \geq 1$
\begin{equation}\label{minsize1}
|B_k(v_i)| \geq |B_1(v_i)| \geq \frac12\log^4m \gg \log^3 n\;.
\end{equation}
Suppose now that $|B_k(v_i)| \leq m/2^{2^d}$. Since we assume that $G$ is a $(\delta,n)$-expander we deduce that
\begin{eqnarray*}
|N_{G'}(B_{k}(v_i))| &\geq& |N_{G}(B_{k}(v_i))| -|W| \geq \frac{\delta2^d }{\log n}|B_{k}(v_i)| -\log^2n\\
&\geq& \frac{\delta2^d }{2\log n}|B_{k}(v_i)|\;,
\end{eqnarray*}
where the last inequality follows from (\ref{minsize1}) and the assumptions that $n \geq n_0(\delta,t)$. We thus get that $G'$ satisfies
(\ref{eq:findpath1}) with $U=\{v_i\}$ and $V=\{v_j\}$ so we can indeed find a path connecting $v_i$ to $v_j$ of length at most $\frac{20}{\delta}\log n \log\log n$. $\qed$

\medskip

\paragraph{Proof of Claim \ref{clm:case21}:}
Let $S_1,\ldots,S_t$ be the sets satisfying the second condition of Claim \ref{clm:sets11}.
Recall that in this case each set $S_i$ is a $\gamma$-expanding ball around a vertex $v_i$ in the induced subgraph $G[S_i]$. Also recall that $m^{1/5} \leq |S_i| \leq m^{1/4}$
and that all the vertices of $G[S_i]$ have degree at most $\log^4 m$.
This means that there is some $k_i$, such that the ball around $v_i$ in the graph $G[S_i]$ satisfies $\log^{4} m \leq |B_{k_i}(v_i)| \leq \log^8 m$.
For every $1 \leq i \leq t$ set $C_i=B_{k_i}(v_i)$ and note that $C_i$ is also a $\gamma$-expanding ball in $G[S_i]$ around the same center $v_i$.

We will shortly show that one can find ${t \choose 2}$ internally vertex disjoint paths $P_{i,j}$,
where each $P_{i,j}$ connects $C_i$ to $C_j$, has length at most
$\frac{20}{\delta}\log n \log\log n$ and avoids all the sets $C_\ell, \ell \not =i,j$.
But let us first observe why this will conclude the proof. For each $1 \leq i \neq j \leq t$ let $Q_{i,j}$ be some path connecting
$v_i$ to the unique vertex of $P_{i,j}$ that belongs to $C_i$. For each $1 \leq i \leq t$ set $V_{i}=\cup_{j}Q_{i,j}$. Then each set $V_i$ is connected in $G$
and the paths $P_{i,j}$ are internally vertex disjoint and avoid the sets $V_\ell, \ell \not =i,j$ so contracting the sets $V_i$
indeed gives us a $K_t$-minor in $G$. As to the order of this minor, note that since $C_i$ is a $\gamma$-expanding-ball (with $\gamma=\delta\log m/5\log n$) of size at most $\log^8 m$
its radius is bounded by
$$
\log_{1+\gamma} |C_i| \leq \frac{2\log(\log^8 m)}{\gamma} \leq \frac{80\log n \log\log m}{\delta\log m}\;,
$$
so each of the paths $Q_{i,j}$ is of length $o(\log n)$ (here we rely on the lemma's assumption that $\log^{1-o(1)} n \leq \log m \leq \log n $).
Hence the total size of the sets $V_i$ is smaller than $\log n$.
Thus, once we show how to find  the above mentioned paths
$P_{i,j}$ with lengths bounded by $\frac{20}{\delta}\log n \log\log n$ we obtain a $K_t$-minor of order at most
$$ \log n+\frac{20t^2}{\delta}\log n\log\log n=O\Big(\frac{t^2}{\delta}\log n\log\log n\Big)\, .$$

Our plan is to show that one can simply iteratively pick the paths $P_{i,j}$ by successive applications of Claim \ref{clm:findpath1}.
To this end, we need to show that after picking some of the paths, we can still find in the remaining graph another path. Since each of the paths
is of length at most $\frac{20}{\delta}\log n \log\log n \ll \log^2 n$, it is enough to show that for any set of vertices $W$ of size at most $\log^2n$, that is disjoint from each of the
sets $C_1,\ldots,C_t$, and such that for any $i \neq j$, the sets $C_i$ and $C_j$ satisfy (\ref{eq:findpath1}) with respect to the graph
\begin{equation}\label{newG1}
G'=G \setminus (W \cup \bigcup_{\ell \neq i,j} C_{\ell})\;.
\end{equation}
Once we know this, we will iteratively pick the paths $P_{i,j}$, where at each iteration we will take $W$ to be the union of the internal vertices of the
paths we have already picked. To establish (\ref{eq:findpath1}) we will need to consider two ``regimes'' of growth; the first is when $C_i$ grows within the set $S_i$, and the second, when it grows out of $S_i$.

Let us first make a simple (but crucial) observation about the sets $C_i$. Recall that $S_i=B_{k'_i}(v_i)$ in $G[S_i]$ for some vertex $v_i$ and $k'_i$,
that $|S_i| \geq m^{1/5}$ and that each vertex in $S_i$ has degree at most $\log^4 m$. This means that $k'_i \geq \log m/20\log\log m$. We also recall that $|C_i| \leq \log^8 m$ and
since $S_i$ is a $\gamma$-expanding ball with $\gamma \geq \delta\log m/5\log n$ this means that in the induced subgraph $G[S_i]$ we have that $C_i=B_{k_i}(v_i)$ for some
$$
k_i \leq \frac{2\log(\log^8 m)}{\gamma} \leq \frac{80\log n \log\log m}{\delta\log m} \leq (\log\log m)^4\;,
$$
where the second inequality relies on the lemma's assumption that $\log m \geq \log^{1-o(1)} n$.
Combining the above facts about $k_i$ and $k'_i$ we get that setting $r_i=\log m/30\log\log m$ we have in the graph $G[S_i]$
\begin{equation}\label{ballinball1}
B_{r_i}(C_i) \subseteq S_i\;.
\end{equation}

We now turn to show that $C_i$ satisfies (\ref{eq:findpath1}) with respect to $G'$ defined in (\ref{newG1}).
Since $S_i$ is a $\gamma$-expanding ball we get from (\ref{ballinball}) that
$|N_{G[S_i]}(C_i)| \geq \delta|C_i|\log m/5\log n$.
Recall now that $S_i$ is disjoint from $\bigcup_{r \neq i,j} C_r$, and that (\ref{ballinball1}) implies that $N_{G[S_i]}(C_i) \subseteq S_i$. Hence
\begin{eqnarray*}
|N_{G'}(C_i)| &\geq&  |N_{G[S_i]}(C_i)|- |W| \geq \delta|C_i|\log m/5\log n- \log^2n \\ &\geq& \delta|C_i|\log m/10\log n\;,
\end{eqnarray*}
where the last inequality uses the fact that $|C_i| \geq \log^4 m$ and the assumption $\log m \geq \log^{1-o(1)} n$. In other words, we have in $G'$
\begin{equation}\label{expandinSi1}
|B_1(C_i)|\geq \left(1+\frac{\delta\log m}{10 \log n}\right)|C_i|\;.
\end{equation}
Recalling (\ref{ballinball1}) we can continue inductively and get that for all $k \leq r_i$ we have
\begin{eqnarray}
\label{beny2}
|N_{G'}(B_{k}(C_i))| &\geq& |N_{G[S_i]}(B_{k}(C_i))| - |W| \geq \delta|B_{k}(C_i)|\log m/5\log n- \log^2n \nonumber\\
&\geq& \delta|B_{k}(C_i)|\log m/10\log n\;,
\end{eqnarray}
where in the last inequality we again used the facts that $|B_{k}(C_i)| \geq |C_i| \geq \log^4 m$ and that $\log m \geq \log^{1-o(1)} n$.
We get that $C_i$ satisfies (\ref{eq:findpath1}) for all $k \leq r_i$. To see this,
note that by (\ref{ballinball1}) $B_{k}(C_i) \subset B_{r_i}(C_i) \subseteq S_i$ implying that $|B_{k}(C_i)| \ll \sqrt{m}$.
This means
that in (\ref{eq:findpath1}) the relevant $d$ is $\log\log m -1$, implying that we should show that $|N_{G'}(B_{k}(C_i))|\geq \delta|B_{k}(C_i)|\log m/20\log n$
as we indeed derive in (\ref{beny2}).
>From (\ref{beny2}) we also have
\begin{equation}\label{largeball1}
|B_{r_i}(C_i)| \geq |C_i|\left(1+\frac{\delta\log m}{10 \log n}\right)^{r_i} \geq \left(1+\frac{\delta\log m}{10 \log n}\right)^{\frac{\log m}{30\log\log m}} \geq \log^{11}m\;,
\end{equation}
where the last inequality uses the the assumption that $\log m \geq \log^{1-o(1)} n$.

Consider now some $k > r_i$ and suppose $|B_k(C_i)| \leq m/2^{2^{d}}$. Then since $G$ is assumed to be a $(\delta,n)$-expander we have
\begin{eqnarray*}
|N_{G'}(B_{k}(C_i))| &\geq& |N_{G}(B_{k}(C_i))| - |W|-\sum_{\ell}|C_{\ell}| \geq \delta2^d|B_{k}(C_i)|/\log n- \log^9m \\
&\geq& \delta2^d|B_{k}(C_i)|/2\log n\;,
\end{eqnarray*}
where the last inequality follows by combining the assumption that $\log m \geq \log^{1-o(1)} n$ with (\ref{largeball1}) which tells us that for any $k \geq r_i$ we have
$|B_{k}(C_i)| \geq |B_{r_i}(C_i)| \geq \log^{11} m$.
So $C_i$ satisfies (\ref{eq:findpath1}) for all $d$.
$\qed$

\end{document}